\documentclass[11pt,a4paper,oneside,leqno]{article}

\usepackage{enumitem}
\usepackage{amsmath,amssymb,amsthm,amsrefs,stmaryrd}
\usepackage[hidelinks]{hyperref}

\hypersetup{
    pdfauthor={Tom\'a\v{s} Crh\'ak},
    pdftitle={On functors from Giry algebras to convex spaces},
    pdfkeywords={convex space, Giry monad},
}

\setenumerate[1]{label=(\roman*)}

\newtheorem{theorem}{Theorem}
\newtheorem{lemma}{Lemma}
\newcommand{\G}{\mathcal{G}}
\newcommand{\Set}{\textbf{\textup{Set}}}
\newcommand{\Meas}{\textbf{\textup{Meas}}}
\newcommand{\Cvx}{\textbf{\textup{Cvx}}}
\newcommand{\SCvx}{\textbf{\textup{SCvx}}}
\newcommand{\Giry}{\Meas^{\G}}
\newcommand{\eval}[1]{ev_{#1}}
\newcommand{\AND}{\:\&\:}
\newcommand{\I}{\mathbb{I}}
\newcommand{\0}{\mathbf{0}}
\newcommand{\1}{\mathbf{1}}
\newcommand{\2}{\mathbf{2}}
\newcommand{\CF}{\Phi}
\newcommand{\SCF}{\Psi}
\newcommand{\sig}{$\sigma$}
\newcommand{\preimage}[1]{#1^{\shortleftarrow}}
\newcommand{\set}[1]{\left\{ #1 \right\}}
\newcommand{\const}{\overline}
\newcommand{\F}[1]{F\left( #1 \right)}
\newcommand{\ch}[1]{\chi_{\raisebox{-.5ex}{$\scriptstyle #1$}}}

\let\phi=\varphi
\let\epsilon=\varepsilon

\begin{document}

\title{On functors from category of Giry algebras to~category of convex spaces}
\author{Tom\'a\v{s} Crh\'ak}
\date{26 Mar 2018}
\maketitle

\begin{abstract}
In~\cite{sturtz2018factorization-v2} the author asserts that the category of
convex spaces is equivalent to the category of Eilenberg-Moore algebras
over the Giry monad. Some of the statements employed in the proof
of this claim have been refuted in our earlier paper~\cite{crhak2018note}.
Building on the results of that paper we prove that no such equivalence exists
and a parallel statement is proved for the category of super convex spaces.
\end{abstract}


\section{Introduction}


The Giry monad $(\G, \eta, \mu)$ on the category of measurable spaces
$\Meas$ has been introduced in~\cite{giry1982categorical}. The functor $\G$ maps
every measurable space $X$ to the measurable space $\G(X)$ of all probability
measures on $X$. The space $\G(X)$ is endowed with the \sig-algebra induced
by the evaluations
\begin{equation*}
    \eval{K}:\G(X)\to\I
    \text{,}
    \qquad
    K\in\Sigma X
    \text{,}
\end{equation*}
where $\Sigma X$ is the \sig-algebra of $X$, and $\I$ is the unit interval.
Every measurable map $f:X\to Y$ is sent to the pushforward map $\G(f)$
defined by
\begin{equation*}
    \G(f)(\phi) = \phi\circ\preimage{f}
    \qquad
    \forall\phi\in\G(X),
\end{equation*}
where $\preimage{f}$ denotes the inverse image map. The unit $\eta:X\to\G(X)$
of the Giry monad assigns the Dirac measure to every element of $X$,
and the counit $\mu:\G^2(X)\to\G(X)$ is given by the integral
\begin{equation*}
    \mu(\phi)(K) = \int_{\xi\in\G(X)}\xi(K) \,d\phi
    \qquad
    \forall\phi\in\G^2(X),K\in\Sigma X.
\end{equation*}

The Eilenberg-Moore algebras over the Giry monad, called
\textit{Giry algebras} throughout this paper, will be written
as pairs $(X,\kappa)$, where $X$ is a measurable space,
and $\kappa:\G(X)\to X$ is the structure map. The category
of Giry algebras will be denoted by $\Giry$.


For a treatment of convex spaces we refer to~\cite{fritz2015convex-v3}
and especially draw the reader's attention to Lemma~3.8, which is employed
in our proof. Herein we stick to the definition of the convex structure
on a set $A$ by means of a~family of binary operations
\begin{equation*}
    +_r:A\times A\to A,
    \qquad
    r\in[0,1],
\end{equation*}
called the \textit{convex combination operations},
satisfying the following axioms:

\begin{description}
\item[unit law]
    \begin{equation*}
        x +_0 y = x
    \end{equation*}

\item[idempotency]
    \begin{equation*}
        x +_r x = x
    \end{equation*}

\item[parametric commutativity]
    \begin{equation*}
        x +_r y = y +_{1-r} x
    \end{equation*}

\item[deformed parametric associativity]
    \begin{equation*}
        (x +_r y) +_s z = x +_t ( y +_{s/t} z )
        \text{,}
    \end{equation*}
    where $s>0$ and $t = r + s - rs$.
\end{description}
Maps of convex spaces preserving the convex combination operations are
called \textit{affine} and the category of convex spaces with affine
maps as morphisms is denoted by $\Cvx$.


\begin{theorem}\label{theorem for cvx}
There is no covariant fully faithful functor from the category of Giry algebras
to the category of convex spaces.
\end{theorem}

The proof of the theorem is deferred to Section~\ref{proof of theorem}.
An outline of the reasoning reads as follows: we define a functor $\CF$
from the category of Giry algebras to the category of convex spaces and
show that the functor is not full. Next we prove that any fully faithful
functor would have to be naturally isomorphic to $\CF$,
but then it cannot be full---a contradiction.

An immediate corollary is that the category of Giry algebras
is not equivalent to the category of convex spaces, refuting the assertion
of~\cite{sturtz2018factorization-v2}.


\section{Preliminaries}


As a Giry algebra, the unit interval $\I$ is endowed
with the Borel \sig-algebra, and its structure map $E$ is given by
\begin{equation*}
    E(\phi) = \int_{x\in\I} x \,d \phi
    \qquad
    \forall\phi\in\G(\I)
    \text{.}
\end{equation*}

More generally, for a set $M$, let $\I^M$ denote the set of \textit{all}
functions from the set $M$ to the unit interval. As a measurable space,
$\I^M$ is endowed with the \textit{evaluation} \sig-algebra,
i.e., the \sig-algebra induced by the evaluation maps
\begin{equation*}
    \eval{m}:\I^M\to\I
    \text{,}
    \qquad
    \forall m\in M
    \text{,}
\end{equation*}
where $\eval{m}(\alpha) = \alpha(m)$. The structure map $E_M$ defined by
\begin{equation*}
    E_M(\phi)(m) = \int_{\I^M} \eval{m} \,d \phi
    \qquad
    \forall\phi\in\G\left(\I^M\right)
    \text{,}
\end{equation*}
makes $(\I^M,E_M)$ a Giry algebra.

\begin{lemma}\label{not measurable}
Let $\const{0}:M\to\I$ be the zero constant function. If the set $M$
is uncountable, then $\set{\const{0}}$ is not measurable in $\I^M$.

\begin{proof}
Let $A$ be the free convex space over $M$. Recall that
$A$ is the set of maps $x\in\I^M$ with $x^{-1}(0)$ cofinite in $M$
and $\sum_{m\in M}x(m) = 1$, and $A$ has the pointwise convex
structure. Consider the set $\Cvx(A,\I)$ of affine maps,
with the evaluation \sig-algebra. Then the map
\begin{equation*}
    f:\I^M\to\Cvx(A,\I)
\end{equation*}
given by
\begin{equation*}
    f(\alpha)(x) = \sum_{m\in M}\alpha(m)x(m)
    \qquad
    \forall\alpha\in\I^M,x\in A
\end{equation*}
yields an isomorphism of measurable spaces, and $f(\const{0})$
is the constant zero function $A\to\I$. Now the conclusion
follows from~\cite{crhak2018note}*{Lemma 4}.
\end{proof}
\end{lemma}


Viewing the unit interval as a convex space, it is endowed with the usual convex
structure
\begin{equation*}
    x +_r y = (1-r)x + ry
    \qquad
    \forall r\in[0,1]
    \quad
    \forall x,y\in\I
    \text{.}
\end{equation*}


Recall that for every pair of elements $x,y$ of a convex space $A$,
the \textit{path map} $\pi_{x,y}:\I\to A$ is defined by
$\pi_{x,y}(r) = x +_r y$, and it is the unique affine map satisfying
\begin{equation*}
    \pi_{x,y}(0) = x
    \AND
    \pi_{x,y}(1) = y
    \text{.}
\end{equation*}


Let us now turn the attention to the two point object $\2=\set{0,1}$. Its
convex structure is determined by
\begin{equation}\label{nat cvx on 2}
    0 +_r 1
        = \begin{cases}
            0 & \text{if $r < 1$,} \\
            1 & \text{otherwise.}
        \end{cases}
\end{equation}
This convex space allows for the definition of the characteristic
function, which we denote by $\ch{}$. In what follows we are going to
characterize the sets whose characteristic function is affine.


A subset of a convex space $A$ is said to be \textit{convex} if it is closed
under the convex combination operations. The empty set, all singletons,
$A$ itself, and the image and the inverse image of a convex set
under an affine map are all convex.

\begin{lemma}\label{affinity of ch}
Let $K$ be a subset of a convex space $A$. Then the characteristic function
$\ch{K}:A\to\2$ is affine if and only if $K$ is convex and satisfies the
following coconvexity condition:
\begin{equation}\label{cocvx}
    x +_r y \in K \implies x\in K
    \qquad
    \forall r < 1
    \quad
    \forall x,y\in A
    \text{.}
\end{equation}
\begin{proof}
First assume that $\ch{K}$ is affine. Then $K$ is the inverse image
of the set $\set{1}\subseteq\2$, hence convex. For $r<1$ and
$x,y\in A$ with $x +_r y \in K$ we have
\begin{equation*}
    \ch{K}(x) +_r \ch{K}(y)
        = \ch{K}(x +_r y)
        = 1
    \text{,}
\end{equation*}
which is possible only if $\ch{K}(x)=1$, thus $x\in K$.

Let us on the other hand suppose that $K$ is convex
and satisfies condition~\eqref{cocvx}. We have the following cases:
\begin{enumerate}
    \item $x,y\in K$.
    Then also $x +_r y \in K$, and
    \begin{equation*}
        \ch{K}(x +_r y)
            = 1
            = 1 +_r 1
            = \ch{K}(x) +_r \ch{K}(y)
        \text{.}
    \end{equation*}

    \item $x,y\notin K$.
    From~\eqref{cocvx} it follows that the complement of $K$ is convex, so that
    $x +_r y \notin K$, and
    \begin{equation*}
        \ch{K}(x +_r y)
            = 0
            = 0 +_r 0
            = \ch{K}(x) +_r \ch{K}(y)
        \text{.}
    \end{equation*}

    \item $x\notin K, y\in K$.
    For $r=1$ we have
    \begin{equation*}
        \ch{K}(x +_r y)
            = \ch{K}(y)
            = \ch{K}(x) +_r \ch{K}(y)
        \text{,}
    \end{equation*}
    whereas for $r<1$ it follows from~\eqref{cocvx} that $x +_r y \notin K$,
    thus
    \begin{equation*}
        \ch{K}(x +_r y)
            = 0
            = 0 +_r 1
            = \ch{K}(x) +_r \ch{K}(y)
        \text{.}
    \end{equation*}

    \item $x\in K, y\notin K$.
    Then employing the previous case we have
    \begin{equation*}
    \begin{split}
        \ch{K}(x +_r y)
            &= \ch{K}(y +_{1-r} x) \\
            &= \ch{K}(y) +_{1-r} \ch{K}(x)
            = \ch{K}(x) +_r \ch{K}(y)
        \text{.}
    \end{split}
    \end{equation*}
\end{enumerate}
\end{proof}
\end{lemma}


As a measurable space, $\2$ is endowed with the discrete \sig-algebra.
The function $h:\G(\2)\to\I$ defined by
\begin{equation*}
    h(\phi) = \phi(\set{1})
    \qquad
    \forall\phi\in\G(\2)
\end{equation*}
is an isomorphism of the Giry algebras $(\G(\2),\mu_\2)$ and $(\I,E)$.
Moreover, the path map $\pi_{0,1}:\I\to\2$ is measurable, and the composition
\begin{equation}\label{def of epsilon_2}
    \epsilon_\2 = \pi_{0,1}\circ h
\end{equation}
provides a structure map, making $(\2, \epsilon_\2)$ a Giry algebra.


\section{Convex structure induced by structure map}\label{induced cvx}


Let $X$ be a measurable space. The underlying set of $\G(X)$ has a natural
\textit{pointwise convex structure}
\begin{equation}\label{nat cvx on G}
    (\phi +_r \psi)(K) = \phi(K) +_r \psi(K)
    \qquad
    \forall K\in\Sigma X
    \text{.}
\end{equation}
Given any measurable map $f:X\to Y$, a simple computation shows that $\G(f)$
is affine with respect to the pointwise convex structures, defined as above,
on $\G(X)$ and $\G(Y)$.


For a Giry algebra $(X,\kappa)$,
the \textit{convex structure induced on $X$ by $\kappa$} is defined by
\begin{equation}\label{cvx for Giry}
    \alpha +_r \beta
    =
    \kappa\bigl( \eta_X(\alpha) +_r \eta_X(\beta) \bigr)
    \text{,}
\end{equation}
for all $r\in[0,1]$ and $\alpha,\beta\in X$.

It is obvious that the operations given by~\eqref{cvx for Giry} satisfy
the unit law, idempotency, and parametric commutativity axioms. The deformed
parametric associativity axiom follows from
\begin{equation*}
    \kappa(\phi +_r \psi) = \kappa(\phi) +_r \kappa(\psi)
    \qquad
    \forall\phi,\psi\in\G(X)
    \text{,}
\end{equation*}
which is in turn an easy consequence of the equality
\begin{equation}\label{mu equality}
    \phi +_r \psi
    =
    \mu_X\bigl( \eta_X(\phi) +_r \eta_X(\psi) \bigr)
    \text{,}
\end{equation}
%
%
where the (pointwise) convex combination operation on the right
is carried out in $\G^2(X)$.


For all morphisms $f:(X,\kappa)\to(Y,\lambda)$ we furthermore have
\begin{align*}
    f(\alpha +_r \beta)
        &= f\bigl( \kappa\bigl( \eta_X(\alpha) +_r \eta_X(\beta) \bigr) \bigr) \\
        &= \lambda\bigl( \G(f)\bigl( \eta_X(\alpha) +_r \eta_X(\beta) \bigr) \bigr) \\
        &= \lambda\bigl(
                \G(f)\bigl( \eta_X(\alpha) \bigr)
                +_r
                \G(f)\bigl( \eta_X(\alpha) \bigr)
            \bigr) \\
        &= \lambda\bigl( \eta_Y\bigl(f(\alpha)\bigr) +_r \eta_Y\bigl(f(\beta\bigr) \bigr) \\
        &= f(\alpha) +_r f(\beta)
        \text{,}
\end{align*}
so that $f$ is affine with respect to the convex structures induced on $X$ and $Y$
by the structure maps $\kappa$ and $\lambda$, respectively.


As a result of the considerations above, it is possible to define
a covariant functor
\begin{equation*}
    \CF:\Giry\to\Cvx
\end{equation*}
so that $\CF(X,\kappa)$ is the set $X$ with the convex structure
induced by $\kappa$, and $\CF(f)=f$ for all morphisms of Giry algebras.


\begin{lemma}\label{CF for I^M and 2}
The functor $\CF$ enjoys the following properties:
\begin{enumerate}
\item $\CF(\I^M, E_M)$ has the pointwise convex structure;
\item $\CF(\2,\epsilon_\2)$ has the usual convex structure
      on $\2$, given by equality~\eqref{nat cvx on 2}.
\end{enumerate}

\begin{proof}
Both assertions follow from straightforward computations:
\begin{enumerate}
\item For all $\alpha,\beta\in\I^M$ and $m\in M$ we have
    \begin{align*}
        (\alpha +_r \beta)(m)
            &= E_M\bigl( \eta_{\I^M}(\alpha) +_r \eta_{\I^M}(\beta) \bigr)(m) \\
            &= \int_{\I^M} \eval{m}
               \,d \bigl( \eta_{\I^M}(\alpha) +_r \eta_{\I^M}(\beta) \bigr)   \\
            &= \left( \int_{\I^M} \eval{m} \,d \eta_{\I^M}(\alpha) \right)
               +_r
               \left( \int_{\I^M} \eval{m} \,d \eta_{\I^M}(\beta) \right)   \\
            &= \alpha(m) +_r \beta(m)
        \text{.}
    \end{align*}
\item With $h$ as in the definition~\eqref{def of epsilon_2} of the map
    $\epsilon_\2$ we have
    \begin{equation*}
        0 +_r 1 = \pi_{0,1}\bigl(
                        h\bigl(\eta_\2(0)\bigr) +_r h\bigl(\eta_\2(1)\bigr)
                  \bigr)
                = \pi_{0,1}(0 +_r 1) = \pi_{0,1}(r)
        \text{.}
    \end{equation*}
\end{enumerate}
\end{proof}
\end{lemma}


\begin{lemma}\label{CF is not full}
The functor $\CF$ is not full.
\begin{proof}
Let $M$ be an arbitrary uncountable set. In view of Lemma~\ref{CF for I^M and 2}
it is possible to write the characteristic map of $\set{\const{0}}$ as
\begin{equation*}
    \ch{\set{\const{0}}}:\CF(\I^M,E_M)\to\CF(\2,\epsilon_\2)
    \text{,}
\end{equation*}
where $\const{0}$ is the constant zero function. The set $\set{\const{0}}$
is convex and a simple argument shows that it satisfies the coconvexity
condition~\eqref{cocvx}, thus $\ch{\set{\const{0}}}$ is affine due to
Lemma~\ref{affinity of ch}.

From Lemma~\ref{not measurable} it however results that $\ch{\set{\const{0}}}$
is not measurable.
\end{proof}
\end{lemma}


\begin{lemma}\label{path maps are Giry morphisms}
Let $(X,\kappa)$ be an arbitrary Giry algebra.
Then all path maps $\pi_{\alpha,\beta}:\I\to\CF(X,\kappa)$
belong to $\Giry\bigl((\I,E),(X,\kappa)\bigr)$.

\begin{proof}
Note that from equality~\eqref{mu equality} it follows that the convex structure
induced on $\G(X)$ by $\mu_X$ is simply the pointwise convex structure
defined by~\eqref{nat cvx on G}, so that for all $\phi,\psi\in\G(X)$ the
path map $\pi_{\phi,\psi}$ is measurable, and
\begin{align*}
    \mu_X\bigl( \G(\pi_{\phi,\psi})(\omega) \bigr)(K)
        &= \int_{\xi\in\G(X)} \xi(K) \,d \G(\pi_{\phi,\psi})(\omega) \\
        &= \int_{r\in\I} \pi_{\phi,\psi}(r)(K) \,d\omega \\
        &= \int_{r\in\I} \bigl( \phi(K) +_r \psi(K) \bigr) \,d\omega \\
        &= \left( \phi(K)\int_{r\in\I} (1-r) \,d\omega \right)
           +
           \left( \psi(K)\int_{r\in\I} r \,d\omega \right) \\
        &= \phi(K) +_{E(\omega)} \psi(K) \\
        &= \pi_{\phi,\psi}\bigl( E(\omega) \bigr)(K),
\end{align*}
thus $\pi_{\phi,\psi}\in\Giry\bigl((\I,E),(\G(X),\mu_X)\bigr)$.
Since
\begin{equation*}
    \pi_{\alpha,\beta}
    =
    \kappa\circ\pi_{\eta_X(\alpha),\eta_X(\beta)},
\end{equation*}
we conclude that $\pi_{\alpha,\beta} \in \Giry\bigl((\I,E),(X,\kappa)\bigr)$.
\end{proof}
\end{lemma}


\section{Proof of Theorem~\ref{theorem for cvx}}\label{proof of theorem}

Let $F:\Giry\to\Cvx$ be a covariant functor, and assume by contradiction that
$F$ is full and faithful.


The set $\Giry(\1,\0)$, where $\0$ and $\1$ are initial and terminal
objects, respectively, is empty. Hence $\Cvx\bigl(\F{\1},\F{\0}\bigr)$
is also empty and it follows that $\F{\0}$ is an initial convex algebra
and that $\F{\1}\ne\emptyset$. For arbitrary $x,y\in\F{\1}$, let $\const{x}$
and $\const{y}$ denote the constant maps $\F{\1}\to\F{\1}$ with the values
$x$ and $y$, respectively. Since the functor $F$ is full, there are morphisms
$u,v\in\Giry(\1,\1)$ such that $\F{u}=\const{x}$ and $\F{v}=\const{y}$,
but as $\1$ is a terminal Giry algebra, $u$ must be equal to $v$,
thus $\const{x}=\const{y}$, and $x=y$. This shows that $\F{\1}$ is a singleton,
and thus a terminal object in $\Cvx$.

The unique element of $\F{\1}$ will be denoted by $*$ in the sequel.
For every Giry algebra $(X,\kappa)$ we now define a map
$\tau_{X,\kappa}:\CF(X,\kappa)\to\F{X,\kappa}$ by
\begin{equation*}
    \tau_{X,\kappa}(\alpha) = \F{\const{\alpha}}(*)\qquad\forall\alpha\in X
    \text{,}
\end{equation*}
where $\const{\alpha}:\1\to X$ denotes the constant map with value $\alpha$.
We claim that $\tau$ is a natural isomorphism, but then, as a consequence
of Lemma~\ref{CF is not full}, the functor $F$ is not full, a contradiction.


To prove the injectivity of $\tau_{X,\kappa}$, consider arbitrary
$\alpha,\beta\in X$. From $\tau_{X,\kappa}(\alpha) = \tau_{X,\kappa}(\beta)$
it follows that $\F{\const{\alpha}} = \F{\const{\beta}}$,
but then $\alpha = \beta$ as $F$ is faithful.
In order to prove that $\tau_{X,\kappa}$ is surjective, take arbitrary
$x\in\F{X,\kappa}$. Since $F$ is full, there exists $\alpha\in X$ such that
$\F{\const{\alpha}}(*) = x$, i.e., $\tau_{X,\kappa}(\alpha) = x$.

Next we show that $\tau$ is natural. Let $u:(X,\kappa)\to(Y,\lambda)$ be
a morphism of Giry algebras. Then for all $\alpha\in X$ we have
\begin{equation*}
\begin{split}
    \F{u}\bigl( \tau_{X,\kappa}(\alpha) \bigr)
        &= \F{u}\bigl( \F{\const{\alpha}}(*) \bigr) \\
        &= \F{u\circ\const{\alpha}}(*)
        = \F{\const{u(\alpha)}}(*)
        = \tau_{Y,\lambda}\bigl( u(\alpha) \bigr).
\end{split}
\end{equation*}


The crux of the proof thus depends in showing that all $\tau_{X,\kappa}$
are affine. Since $\tau_{\I,E}$ is bijective, we may, for all $r\in[0,1]$,
define operations $\oplus_r$ by the formula
\begin{equation*}
    x \oplus_r y
        = \tau_{\I,E}^{-1} \bigl( \tau_{\I,E}(x) +_r \tau_{\I,E}(y) \bigr)
    \qquad
    \forall x,y\in\I
    \text{.}
\end{equation*}
Straightforward calculations show that this family of operations defines
a convex structure on $\I$. Moreover, for all $u\in\Giry(\I,\I)$ we have
\begin{equation*}
\begin{split}
    \tau_{\I,E}\bigl( u(x \oplus_r y) \bigr)
        &= (\tau_{\I,E}\circ u)\bigl( \tau_{\I,E}^{-1}
                \bigl( \tau_{\I,E}(x) +_r \tau_{\I,E}(y) \bigr)
            \bigr)                                               \\
        &= \F{u}\bigl( \tau_{\I,E}(x) +_r \tau_{\I,E}(y) \bigr)   \\
        &= \F{u}\bigl( \tau_{\I,E}(x) \bigr)
           +_r
           \F{u}\bigl( \tau_{\I,E}(y) \bigr)                      \\
        &= \tau_{\I,E}\bigl( u(x) \bigr)
           +_r
           \tau_{\I,E}\bigl( u(y) \bigr)
    \text{,}
\end{split}
\end{equation*}
hence, since $\tau_{\I,E}$ is injective,
\begin{equation*}
    u(x \oplus_r y) = u(x) \oplus_r u(y)
    \text{.}
\end{equation*}
From Lemma~\ref{path maps are Giry morphisms} it follows that for all
$a,b\in\I$ the path maps $\pi_{a,b}$ are elements of
$\Giry(\I,\I)$, thus
\begin{equation*}
    \pi_{a,b}(x \oplus_r y) = \pi_{a,b}(x) \oplus_r \pi_{a,b}(y)
    \text{,}
\end{equation*}
and we apply~\cite{fritz2015convex-v3}*{Lemma 3.8} to conclude
that $\oplus_r$ and $+_r$ in fact coincide, but then
\begin{equation*}
    \tau_{\I,E}(x +_r y) = \tau_{\I,E}(x) +_r \tau_{\I,E}(y)
    \text{.}
\end{equation*}
For every Giry algebra $(X,\kappa)$ and $\alpha,\beta\in X$ we therefore have
\begin{equation*}
\begin{split}
    \tau_{X,\kappa}(\alpha +_r \beta)
        &= \tau_{X,\kappa}\bigl( \pi_{\alpha,\beta}(r) \bigr)                     \\
        &= \F{\pi_{\alpha,\beta}}\bigl( \tau_{\I,E}(r) \bigr)                     \\
        &= \F{\pi_{\alpha,\beta}}\bigl( \tau_{\I,E}( 0 +_r 1 ) \bigr)             \\
        &= \F{\pi_{\alpha,\beta}}\bigl( \tau_{\I,E}(0) +_r \tau_{\I,E}(1) \bigr)  \\
        &= \F{\pi_{\alpha,\beta}}\bigl( \tau_{\I,E}(0) \bigr)
           +_r
           \F{\pi_{\alpha,\beta}}\bigl( \tau_{\I,E}(1) \bigr)                     \\
        &= \tau_{X,\kappa}\bigl( \pi_{\alpha,\beta}(0) \bigr)
           +_r
           \tau_{X,\kappa}\bigl( \pi_{\alpha,\beta}(1) \bigr)                     \\
        &= \tau_{X,\kappa}(\alpha) +_r \tau_{X,\kappa}(\beta)
    \text{,}
\end{split}
\end{equation*}
which concludes the proof.


\section[Appendix: Super convex spaces]%
{Appendix: Super convex spaces\footnote{Added in Sep 2020}}

Similar statement can be made for super convex spaces---the proof reads along
the same lines as the one given above, we however need to show that the natural
isomorphism $\tau$ is countably affine.


In what follows, $\omega$ denotes the set of natural numbers and
$\Omega\subseteq\I^\omega$
is the set of all sequences summing to one. In particular, let us,
for all $n\in\omega$, define elements $\delta^n\in\Omega$ to be the
unique sequences of $\Omega$ with $\delta^n_n = 1$.

A \textit{super convex space} $A$ is a set $A$ together with a structure
\begin{equation*}
    \Omega\times A^\omega\to A
    \text{,}
\end{equation*}
which maps all $r\in\Omega$ and $a\in A^\omega$ to elements denoted by
$\sum_n r_n a_n$, satisfying
\begin{enumerate}
    \item $\sum_m \delta^n_m a_m = a_n$,
    \item $\sum_n r_n \bigl( \sum_m s^n_m a_m \bigr)
            = \sum_m \bigl( \sum_n r_n s^n_m \bigr) a_m$.
\end{enumerate}

Given super convex spaces $A$, $B$, a map $f\in\Set(A,B)$ is
\textit{countably affine} provided it preserves the super
convex structures of $A$ and $B$, that is
\begin{equation*}
    f\Bigl( \sum r_n a_n \Bigr) = \sum r_n f(a_n)
    \text{.}
\end{equation*}

Let $\SCvx$ denote the category of super convex spaces with the countably
affine maps as morphisms.


Every Giry algebra $(X,\kappa)$ has a natural super convex structure given by
\begin{equation*}
    \sum_n r_n a_n = \kappa\left( \sum_n r_n \eta_X(a_n) \right)
    \qquad
    \forall a\in X^\omega
    \text{,}
\end{equation*}
which we employ to define a covariant functor
\begin{equation*}
    \SCF:\Giry\to\SCvx
    \text{.}
\end{equation*}

On the other hand, each super convex space has the obvious convex structure.
For each Giry algebra $(X,\kappa)$ this convex structure on $\SCF(X,\kappa)$
coincides with the convex structure induced on $X$ by $\kappa$, as defined
in~Section~\ref{induced cvx}.


Observe that $\Omega$ itself is a super convex space, the structure being
defined pointwise. If we consider $\omega$ with the discrete \sig-algebra,
$\SCF(\G(\omega),\mu_\omega)$ and $\Omega$
are isomorphic in $\SCvx$; explicitly, the isomorphism
\begin{equation*}
    d:\Omega\to\G(\omega)
\end{equation*}
is given by
\begin{equation*}
    d(r)(K) = \sum_{n\in K}r_n
    \qquad
    \forall r\in\Omega,K\subseteq\omega
    \text{.}
\end{equation*}
Moreover, when we endow $\Omega$ with the initial \sig-algebra generated
by the evaluations
$\eval{n}:\Omega\to\I$, $n\in\omega$,
then $d$ becomes an isomorphism of Giry algebras.


Let $A$ be a super convex space. In analogy to the path maps
$\pi_{x,y}$,
let us define, for all $a\in A^\omega$,
the map $\pi_a:\Omega\to A$ by
\begin{equation*}
    \pi_a(r) = \sum r_n a_n
    \qquad
    \forall r\in\Omega
    \text{,}
\end{equation*}
which is in fact countably affine. Given an arbitrary Giry algebra $(X,\kappa)$,
it is easy to check that for $A=\SCF(X,\kappa)$ we have
\begin{equation*}
    \pi_a = \kappa\circ\G(a)\circ d
    \text{,}
\end{equation*}
and thus $\pi_a\in\Giry(\Omega,(X,\kappa))$.


\begin{lemma}\label{SCF is not full}
The functor $\SCF$ is not full.
\begin{proof}
The characteristic map $\ch{\set{\const{0}}}$ considered
in the proof of~Lemma~\ref{CF is not full} is countably affine.
\end{proof}
\end{lemma}


\begin{lemma}\label{characterization of countably affine maps}
Let $A$ and $B$ be super convex spaces and $f\in\Set(A,B)$. Then
$f$ is countably affine if and only if $f\circ\pi_a$ is countably affine
for all $a\in A^\omega$.

\begin{proof}
For all $r\in\Omega$ and $a\in A^\omega$ we have
\begin{equation*}
    r = \sum r_n \delta^n
    \quad
    \text{and}
    \quad
    a_n = \pi_a(\delta^n)
    \text{,}
\end{equation*}
therefore if $f\circ\pi_a$ is countably affine for all $a\in A^\omega$,
then
\begin{equation*}
\begin{split}
    f\Bigl( \sum r_n a_n \Bigr)
        &= f\bigl( \pi_a(r) \bigr) \\
        &= (f\circ\pi_a)\Bigl( \sum r_n \delta^n \Bigr)
         = \sum r_n f\bigl( \pi_a(\delta^n) \bigr)
         = \sum r_n f(a_n)
    \text{.}
\end{split}
\end{equation*}
\end{proof}
\end{lemma}


\begin{lemma}\label{countably affine on I}
Let $A$ be an arbitrary super convex space and $f\in\Set(\I,A)$. Then
$f$ is countably affine if and only if $f$ is affine, i.e.,
\begin{equation*}
    \SCvx(\I,A) = \Cvx(\I,A)
    \text{.}
\end{equation*}

\begin{proof}
For $r\in\Omega$ and $a\in\I^\omega$ let $x = \sum r_n a_n$. Then we have
\begin{equation*}
\begin{split}
    f\Bigl( \sum r_n a_n \Bigr)
        &= f(0 +_x 1)
         = f(0) +_x f(1) \\
        &= \sum r_n \bigl( f(0) +_{a_n} f(1) \bigr)
         = \sum r_n f(0 +_{a_n} 1)
         = \sum r_n f(a_n)
    \text{.}
\end{split}
\end{equation*}
\end{proof}
\end{lemma}


\begin{theorem}
There is no covariant fully faithful functor from the category of Giry algebras
to the category of super convex spaces.

\begin{proof}
Let $F:\Giry\to\SCvx$ be a functor and let us assume by contradiction that
$F$ is full and faithful. As in the proof for convex spaces
in~Section~\ref{proof of theorem}, one constructs natural maps
$\tau_{X,\kappa}:\SCF(X,\kappa)\to F(X,\kappa)$, which are isomorphisms
of the underlying convex spaces.

We will show that $\tau_{X,\kappa}$ is countably affine for all
Giry algebras $(X,\kappa)$. For $\tau_{\I,E}$ this fact follows
from~Lemma~\ref{countably affine on I}. As for $\tau_\Omega$,
by naturality of $\tau$ we have for all $m\in\omega$
\begin{equation*}
    \tau_{\I,E}\circ\eval{m}
    =
    \F{\eval{m}}\circ\tau_\Omega
    \text{,}
\end{equation*}
so that for all $r\in\Omega$ and $s\in\Omega^\omega$
\begin{equation*}
\begin{split}
    (\tau_{\I,E}\circ\eval{m}) \Bigl( \sum r_n s_n \Bigr)
        &= \sum r_n \tau_{\I,E}\bigl( \eval{m}(s_n) \bigr)                      \\
        &= \sum r_n \F{\eval{m}}\bigl( \tau_\Omega(s_n) \bigr)                  \\
        &= \F{\eval{m}}\Bigl( \sum r_n \tau_\Omega(s_n) \Bigr)                  \\
        &= \bigl( \F{\eval{m}}\circ\tau_\Omega \bigr)
           \biggl(
                \tau_\Omega^{-1} \Bigl( \sum r_n \tau_\Omega(s_n) \Bigr)
           \biggr)                                                              \\
        &= (\tau_{\I,E}\circ\eval{m})
           \biggl(
                    \tau_\Omega^{-1}\Bigl( \sum r_n \tau_\Omega(s_n) \Bigr)
           \biggr)
\end{split}
\end{equation*}
and we employ the injectivity of $\tau_{\I,E}$ to conclude
\begin{equation*}
    \tau_\Omega\Bigl( \sum r_n s_n \Bigr)
    =
    \sum r_n \tau_\Omega(s_n)
    \text{,}
\end{equation*}
that is, $\tau_\Omega$ is countably affine.

In general, let $(X,\kappa)$ be an arbitrary Giry algebra.
By naturality of $\tau$, for all $a\in X^\omega$, we have
\begin{equation*}
    \tau_{X,\kappa}\circ\pi_a
    =
    \F{\pi_a}\circ\tau_\Omega
    \text{,}
\end{equation*}
and therefore
\begin{equation*}
    \tau_{X,\kappa}\circ\pi_a\in\SCvx(\Omega, \F{A})
    \text{,}
\end{equation*}
so that $\tau_{X,\kappa}$ is countably affine
by~Lemma~\ref{characterization of countably affine maps}.

It follows that $\tau$ is a natural isomorphism of $\SCF$ and $F$, but
then $F$ cannot be full due to~Lemma~\ref{SCF is not full}, a contradiction.
\end{proof}
\end{theorem}


\pagebreak

\end{document}